\documentclass{article}

\usepackage[english]{babel}
\usepackage[letterpaper,top=2cm,bottom=2cm,left=3cm,right=3cm,marginparwidth=1.75cm]{geometry}

\usepackage{blindtext}
\usepackage{titlesec}
\usepackage{amsmath}
\usepackage{amssymb}
\usepackage{amsfonts}
\usepackage{graphicx}
\usepackage[colorlinks=true, allcolors=blue]{hyperref}
\usepackage{bbm}
\usepackage{mathabx}
\usepackage{amsthm}
\usepackage{tikz} 
\usepackage{mathtools}
\usepackage{hyperref}
\usepackage[capitalise]{cleveref}
\usepackage{animate}
\usepackage{xcolor}
\usepackage{graphicx}
\usepackage{xparse}
\usepackage{enumitem}
\usepackage{titlesec}
\usepackage{bbm}
\usepackage{tikz} 

\titleformat*{\section}{\LARGE\bfseries}
\titleformat*{\subsection}{\Large\bfseries}

\newtheorem{theorem}{Theorem}
\newtheorem*{theorem*}{Theorem}
\newtheorem{lemma}[theorem]{Lemma}

\theoremstyle{definition}
\newtheorem{definition}[theorem]{Definition}

\theoremstyle{remark}
\newtheorem*{remark*}{Remark}

\DeclareMathOperator{\Var}{Var}
\DeclareMathOperator{\Cov}{Cov}

\title{Microscopic Probabilities for $U$-Statistics and \\ Almost-Balanced Permutations}
\author{Hillel Schein\thanks{Department of Mathematics,
Technion---Israel Institute of Technology, Haifa, Israel.\\
Email: \href{mailto:hillel.raz@campus.technion.ac.il}
{\texttt{hillel.raz@campus.technion.ac.il}}.}}
\date{November 2026}
\hypersetup{pdftitle={Local steering of U-statistics},pdfauthor={Hillel Schein}}

\begin{document}
\maketitle
\begin{abstract}
We estimate the probability that a vector of U-statistics
lies within $C/n$ of a target in the range of typical
fluctuations. For sufficiently large $C$, we prove
a lower bound of the typical Gaussian order under finite variance,
and a matching upper bound for bounded kernels. We show
that under finite variance alone, stronger concentration
may occur, and find the best possible asymptotic upper bound. Our proof is algorithmic, estimating probabilities by counting
successful modifications produced by a random algorithm.

We apply this result to settle a question of Beniamini, Lavee and Linial~\cite{BLL25}, who asked if there exist permutations with $k$-order sub-pattern counts at most $O(n^{k-1})$ away from one another for $k>4$. We prove
not only existence, but enumerate these permutations
up to a constant factor, showing that there are
$\Theta(n!n^{-(k-1)^2/2})$ such permutations of size $n$.

\end{abstract}



\section{Introduction}

Given independent random variables $X_{1},\dots,X_{n}$ with common distribution $\mu$ and a symmetric function $h(x_{1},\dots,x_{k})$, known as a \textit{kernel}, the corresponding normalized $U$-statistic is the average
$$U_{n}(h)=\binom{n}{k}^{-1}\sum_{i_{1}<\dots<i_{k}\in [n]} h(X_{i_{1}},\dots,X_{i_{k}}).$$
$U$-statistics are a well-studied family of statistics, with examples including empirical moments~\cite{GD18,Hef97} and sub-pattern densities in permutations~\cite{JNZ15,EvenZohar20}.

When the kernel $h$ is of finite variance, its $U$-statistic is typically studied through its Hoeffding decomposition~\cite{Hoeffding48}, where different components fluctuate at different orders. In particular, we will utilize the first-order component:
$U_{n}(h)-\mathbbm{E}[U_{n}(h)]=\frac{1}{n}\sum_{i\in[n]}(h)_{1}(X_{i})+R_{n}.$

When nontrivial, this first-order component dominates the fluctuations of the $U$-statistic. The sum $\frac{1}{n}\sum_{i\in[n]}(h)_{1}(X_{i})$ fluctuates at order $\frac{1}{\sqrt{n}}$, while the remainder $R_{n}$ fluctuates at order at most $n^{-1}$; see \cref{sec:prelim} for details. Given a vector of $U$-statistics $U_{n}=(U_{n}(h_{j}))_{j\in[d]}$, we denote its first-order fluctuation space by
$S=\operatorname{range}\Cov(((h_{j})_{1}(X_{1}))_{j\in[d]})$, and its dimension by $s=\dim S$.

Thus, denoting $P_{S},Q_{S}$ as the orthogonal projections onto $S,S^{\bot}$, respectively, we have that the projection $P_{S}(U_{n}-\mathbbm{E}[U_{n}])$ fluctuates at order $\frac{1}{\sqrt{n}}$, while $Q_{S}(U_{n}-\mathbbm{E}[U_{n}])$ fluctuates at order at most $n^{-1}$. Our main result gives bounds on the probability that $U_{n}$ lies in small sets.

\begin{theorem}[Microscopic Probabilities]\label{thm:main}
Let $h_{j}(x_{1},\dots,x_{k_{j}}),j\in[d]$ be fixed symmetric kernels with finite variance, and let $U_{n}=(U_{n}(h_{j}))_{j\in[d]}$ be their corresponding vector of $U$-statistics. Fix $M\geq 0$.

Then there exist constants $C,c>0$ such that, for all sufficiently large $n$ and every $z\in S$ satisfying $\|z\|\leq \frac{M}{\sqrt{n}}$, we have
$$P(\|U_{n}-\mathbbm{E}[U_{n}]-z\|\leq\frac{C}{n})\geq cn^{-\frac{s}{2}}.$$
If the kernels $h_{j},j\in[d]$ are also bounded, then for every fixed $C>0$ an upper bound of the same order holds. Explicitly, there exists $C'>0$ such that, for all sufficiently large $n$,
$$\sup_{\substack{z\in S\\ \|z\|\leq \frac{M}{\sqrt{n}}}}P(\|U_{n}-\mathbbm{E}[U_{n}]-z\|\leq\frac{C}{n})\leq C'n^{-\frac{s}{2}}.$$
When the kernels are only known to have finite variance, then for every fixed $C>0$ we have the following upper bounds, which are weaker when $s>2$:
$$\sup_{\substack{z\in S\\ \|z\|\leq \frac{M}{\sqrt{n}}}}P(\|U_{n}-\mathbbm{E}[U_{n}]-z\|\leq\frac{C}{n})=\begin{cases}
    O(n^{-\frac{s}{2}}) & s\in\{0,1,2\},\\
    o(n^{-1}) & s>2.
\end{cases}$$
Further, when $s>2$, the $o(n^{-1})$ upper bound admits no universal rate of improvement. More precisely, for every integer $s>2$, every fixed $C>0$, and every positive sequence $\epsilon_{n}\rightarrow 0$, there exist fixed square-integrable kernels $h_{1},\dots,h_{d}$ whose first-order fluctuation space has dimension $s$, such that
$$\limsup_{n\rightarrow\infty}\frac{n}{\epsilon_{n}}P(\|U_{n}-\mathbbm{E}[U_{n}]\|\leq\frac{C}{n})=\infty.$$
\end{theorem}

Note that $z\in S$, and we study the event $\|U_{n}-\mathbbm{E}[U_{n}]-z\|\leq \frac{C}{n}$. Therefore, whenever this event holds we must have $\|Q_{S}(U_{n}-\mathbbm{E}[U_{n}])\|\leq\frac{C}{n}$. Thus, the result on microscopic probabilities controls the degenerate subspace as well, allowing it to apply when the statistics have a singular first-order covariance.

One such case is the density of sub-patterns in a uniformly random permutation. Applying \cref{thm:main} in this setting, we settle a question of Beniamini, Lavee and Linial~\cite{BLL25} about balanced permutations. Given a permutation $\pi\in S_{n}$ and an integer $k<n$, we say that $\pi$ is $k$-balanced if the number of copies of every $\tau\in S_{k}$ in $\pi$ is the same. Writing $\#\tau(\pi)$ for the number of copies of $\tau$, this is equivalent to
$$D_{k}(\pi)=\max_{\tau\in S_{k}}\left|\#\tau(\pi)-\frac{1}{k!}\binom{n}{k}\right|=0.$$
They constructed such permutations for $k\leq 3$ under the necessary divisibility conditions. For $k\geq 4$, they proved that no such permutation exists and quantified the minimal gap from balance:
$$\min_{\pi\in S_{n}}D_{k}(\pi)=\Omega_{k}(n^{k-1}),\qquad k\geq 4.$$
They showed this bound is tight for $k=4$, and asked if this holds for $k>4$ as well. In a uniformly random permutation, these sub-pattern densities are $U$-statistics (see \cref{sec:prelim}), and the dimension of the first-order fluctuation space of all patterns of order $k$ is $(k-1)^{2}$~\cite{JNZ15}. Applying our result therefore gives both existence and enumeration of almost-balanced permutations up to constant factors.
\begin{theorem}[Almost-Balanced Permutations]\label{thm:permutations}
For every fixed $k\ge2$, there is $C_k>0$ such that, for every fixed $C\ge C_k$, as $n\rightarrow\infty$,
$$\left|\{\pi\in S_n:D_k(\pi)\le Cn^{k-1}\}\right|
 =\Theta_{k,C}\left(n!\,n^{-(k-1)^2/2}\right).$$
Consequently, for every fixed $k\ge4$,
$$\min_{\pi\in S_n}D_k(\pi)=\Theta_k(n^{k-1}).$$
\end{theorem}

\section{Proof Overview and Discussion}\label{sec:overview}

Our proof of \cref{thm:main} takes an algorithmic approach. In \cref{sec:algorithm}, we design a random algorithm to steer the values of $U_{n}$ of a given sample by selectively resampling variables. We choose a basis of $S$, signifying directions in which we can shift $U_{n}$ at first order, and for each such direction find variables for which resampling induces a first-order change in that direction. Given a random set of $L=\Theta(\sqrt{n})$ variables of each such type, each collection of prefixes $(t_{1},\dots,t_{s})\in\{-L,\dots,L\}^{s}$ defines a set of variables to resample in order to shift $U_{n}$. This gives $\Theta(n^{\frac{s}{2}})$ possible modifications of $U_{n}$.

If all modifications in this family preserved the distribution of the sample, we would have that the probability that $U_{n}$ belongs to a prescribed set of values would be exactly the expected proportion of modifications which shift $U_{n}$ to that set. While our algorithm doesn't preserve the distribution exactly, the fact that we only resample $\Theta(\sqrt{n})$ variables means that we don't change the distribution by much. In \cref{sec:probability-proof}, we use this to bound the original probability from both sides by the expected proportions of successful steering options, with suitable restrictions to typical starting configurations. In \cref{sec:control-change}, we analyze the effect of the algorithm directly, and use this analysis in \cref{sec:main-proof} to bound the expected proportion of successful steering options. We show the bounds are of optimal order in \cref{sec:sharpness} by providing a family of examples. 

Estimating probabilities through local modifications has received much attention in recent years, with several works developing general formulations and emphasizing their broader applicability. A particularly close comparison, which helped inspire our proof method, is the work of Fox, Kwan and Sauermann~\cite{FKS21}. They prove upper bounds for small-interval probabilities of functions of independent Bernoulli variables. Their proof follows a random ordering of bit changes, bounds the number of visits to a small interval, and converts this count into a probability bound. Our signed-prefix construction uses a related counting argument in several directions, while supplying lower bounds using the number of visits to a small set as well. 
Further examples include upper bounds on small-interval probabilities for general statistics, first- and last-passage percolation, and directed polymers~\cite{Cha19,BC20,Elb25}; matching upper and lower bounds for point probabilities of edge counts in random induced subgraphs of Ramsey graphs using switching arguments~\cite{KSSS23} (see also~\cite{HM10}); anti-concentration bounds for functions of independent real variables~\cite{AGKW09}; and variance lower bounds for longest common subsequences~\cite{LM09}.

Taken together, these works illustrate the flexibility
of local-modification arguments for studying small-ball
probabilities and fluctuations. By choosing modifications
adapted to the statistic and estimating probabilities
through the number of successful modifications, this
approach can provide bounds without relying on a local
limit theorem, even when no uniform Gaussian-order description of microscopic probabilities is available. Our work further illustrates this
flexibility for vectors of U-statistics. The signed-prefix
construction yields matching upper and lower bounds for
bounded kernels, while under finite variance alone it
also gives sharp general bounds in a setting where
microscopic probabilities can differ substantially from
the Gaussian order. We believe that this adaptability
makes local-modification methods promising for other
settings in which no universal microscopic behavior
can be expected.

\section{Preliminaries}\label{sec:prelim}

\subsection{First-order decomposition and fluctuation scales}\label{sec:first-order}

We begin by formally defining a $U$-statistic.

\begin{definition}
    Let 
    \begin{itemize}
        \item $(\Omega,\mathcal{F},\mu)$ be any underlying probability space.
        \item $k\in\mathbb{N}$ be the order of the statistic, and $h:\Omega^{k}\rightarrow\mathbb{R}$ be a symmetric measurable function. 
        \item $k\leq n\in\mathbb{N}$, and $X_{1},\dots,X_{n}$ be independent $\Omega$-valued random elements with common law $\mu$.
    \end{itemize} 
    Then the corresponding normalized $U$-statistic is
    $$U_{n}(h)=\binom{n}{k}^{-1}\sum_{I\in \binom{[n]}{k}} h(X_{I}),\qquad X_{I}=(X_{i})_{i\in I}.$$
\end{definition}
Note that the ordering of $X_{I}$ does not matter as $h$ is symmetric.

Whenever $\Var(h)<\infty$, the scaled first order component is defined as
\[
(h)_1(x)=k\left(\mathbbm{E}[h(X_{[k]})\mid X_{1}=x]
-\mathbbm{E}[h(X_{[k]})]\right).
\]
Note that $(h)_{1}$ is centered, and that the definition is invariant to the choice of $k$ variables as they are i.i.d. This gives us by definition that $U_{n}(h)-\mathbbm{E}[U_{n}(h)]=\frac{1}{n}\sum_{i\in[n]}(h)_{1}(X_{i})+U_{n}(g_{h})$, where we define
$$g_{h}(x_1,\dots,x_k)=h(x_1,\dots,x_k)
-\mathbbm{E}[h(X_1,\dots,X_k)]
-\frac{1}{k}\sum_{a=1}^{k}(h)_1(x_a).$$

Therefore the remainder is also a $U$-statistic. In the next lemma, we prove that the remainder is degenerate, and show the orders of fluctuation.
\begin{lemma}[First-order decomposition and fluctuation scales]\label{lem:first-order-scales}
Let $h_{j},j\in[d]$ be fixed kernels of orders $k_{j}$ with finite variance, let $n>\max_{j\in[d]}k_{j}$, and let $U_{n}=(U_{n}(h_{j}))_{j\in[d]}$ be their normalized $U$-statistic vector. Then for all $j\in[d]$, the remainder $g_{h_{j}}$ is degenerate, meaning that for any $I\in\binom{[n]}{k_{j}},i\in[n]$ we have
$$\mathbbm{E}[g_{h_{j}}(X_{I})\mid X_{i}]\equiv 0\hspace{2mm}a.s.$$
Additionally, $U_{n}$ satisfies
\[
U_n-\mathbbm{E}[U_n]
=\frac{1}{n}\sum_{i=1}^{n}
\bigl((h_j)_1(X_i)\bigr)_{j\in[d]}+R_n,
\qquad
\mathbbm{E}\bigl[\|R_n\|^2\bigr]=O(n^{-2}).
\]
Here $R_n=(U_n(g_{h_j}))_{j\in[d]}$. Finally, with $S$, $s$, $P_S$ and $Q_S$ as in the introduction, we have
\[
\mathbbm{E}\left[\bigl\|P_S(U_n-\mathbbm{E}[U_n])\bigr\|^2\right]
=O(n^{-1}),
\qquad
\mathbbm{E}\left[\bigl\|Q_S(U_n-\mathbbm{E}[U_n])\bigr\|^2\right]
=O(n^{-2}).
\]
The implied constants depend only on the fixed kernels and their orders.
\end{lemma}

\begin{proof}
Fix $j\in[d]$ and $I\in\binom{[n]}{k_j}$. If $i\in I$, symmetry and the definition of $(h_j)_1$ give
\[
\mathbbm{E}[h_j(X_I)\mid X_i]
=\mathbbm{E}[h_j(X_I)]+\frac{1}{k_j}(h_j)_1(X_i).
\]
All other first-order summands in $g_{h_j}(X_I)$ have conditional expectation zero by independence and centering. Therefore,
\[
\mathbbm{E}[g_{h_j}(X_I)\mid X_i]
=\mathbbm{E}[h_j(X_I)\mid X_i]
-\mathbbm{E}[h_j(X_I)]-\frac{1}{k_j}(h_j)_1(X_i)=0.
\]
If $i\notin I$, the same conclusion follows because $g_{h_j}(X_I)$ is centered and independent of $X_i$.

This degeneracy is what makes $R_{n}$ have small variance. Note that the remainder is centered, and so it is equivalent to estimate the covariance matrix. By bilinearity we have
$$\Cov(R_n)_{j,j'}=\frac{1}{\binom{n}{k_j}\binom{n}{k_{j'}}}
\sum_{I\in\binom{[n]}{k_j},I'\in\binom{[n]}{k_{j'}}}
\Cov(g_{h_{j}}(X_{I}),g_{h_{j'}}(X_{I'})).$$
When $I\cap I'=\emptyset$, the covariance is zero as different variables are independent. When $I\cap I'=\{i\}$, we have by the law of total covariance that
$$\Cov(g_{h_{j}}(X_{I}),g_{h_{j'}}(X_{I'}))=\mathbbm{E}[\Cov(g_{h_{j}}(X_{I}),g_{h_{j'}}(X_{I'})\mid X_{i})]+\Cov(\mathbbm{E}[g_{h_{j}}(X_{I})\mid X_{i}],\mathbbm{E}[g_{h_{j'}}(X_{I'})\mid X_{i}]).$$
The first summand is zero as they only share $X_{i}$ and are thus independent once conditioned on it, and the second summand is zero by degeneracy. Therefore, we only need to consider pairs intersecting in at least two variables for the variance. Since $g_{h_j}$ and $g_{h_{j'}}$ have finite variance, by Cauchy--Schwarz the covariance of any such pair of summands is $O(1)$, and the proportion of these pairs is $O(n^{-2})$, giving the desired bound.

To estimate the fluctuation of the first order component, note that $((h_j)_1(X_i))_{j\in[d]}$ are centered and independent across $i$, so the second moment of their average is
\[
\frac1n\mathbbm{E}\left[\left\|\bigl((h_j)_1(X_1)\bigr)_{j\in[d]}\right\|^2\right].
\]
Combining this with the bound on $R_n$ proves the estimate for $P_S$. If a vector belongs to $S^\bot$, its inner product with the first-order vector has mean and variance zero. Taking an orthonormal basis of $S^\bot$ therefore gives $Q_S\bigl((h_j)_1(X_i)\bigr)_{j\in[d]}=0$ almost surely. Thus $Q_S(U_n-\mathbbm{E}[U_n])=Q_SR_n$, proving the last bound, including when $s=0$.
\end{proof}

\subsection{Permutation patterns as $U$-statistics}\label{sec:permutation-prelim}

We use the following representation of a uniformly random permutation of order $n$, useful for studying sub-pattern densities as $U$-statistics and used in both~\cite{JNZ15,EvenZohar20}. Let $(X_i,Y_i)$, $i\in[n]$, be independent uniform points of $[0,1]^2$. Order them by their first coordinates and read the relative ranks of their second coordinates to obtain the random permutation $\Pi_n$. Explicitly, if $X_{i_1}<\dots<X_{i_n}$, define
$$\forall j\neq k\in[n],\qquad\Pi_{n}(j)<\Pi_{n}(k)\iff Y_{i_{j}}<Y_{i_{k}}.$$
Ties have probability zero, and so this is well defined. The first coordinates are independent of the second coordinates, whose ranks are uniformly distributed over $S_n$; hence $\Pi_n$ is uniform.

For $\tau\in S_k$, let $h_\tau$ be the indicator that $k$ points, ordered by their first coordinates, have second-coordinate pattern $\tau$. This defines a bounded symmetric measurable kernel. The $k$-subsets of the point sample correspond to the $k$-subsets of positions in $\Pi_n$, preserving their induced patterns. Thus
\[
U_n(h_\tau)=\frac{\#\tau(\Pi_n)}{\binom{n}{k}},
\qquad
\mathbbm{E}[U_n(h_\tau)]=\frac{1}{k!}.
\]
To apply \cref{thm:main} to $\Pi_{n}$, we need to know $s$, the dimension of the first-order fluctuation space for $k$-order sub-pattern densities. This was first proved by Janson, Nakamura, and Zeilberger~\cite{JNZ15}, and the finer decomposition of the pattern profile was then characterized by Even-Zohar~\cite{EvenZohar20}. We include an elementary verification of the first-order space dimension for completeness.

\begin{lemma}[First-order rank of the permutation pattern profile]\label{lem:permutation-rank}
For every fixed $k\geq2$, the vector $(U_n(h_\tau))_{\tau\in S_k}$ has first-order fluctuation rank
\[
s=(k-1)^2.
\]
\end{lemma}

\begin{proof}
For $i\in[k]$, define $b_{i}(t)=\binom{k-1}{i-1}t^{i-1}(1-t)^{k-i}$. Conditional on one point being $(x,y)$, its ranks out of any given $k$ points in the two coordinates are independent, with probabilities $b_{i}(x)$ and $b_{j}(y)$. Given ranks $i,j$, the correspondence between the remaining first- and second-coordinate ranks is a uniform bijection. It produces $\tau$ with probability $1/(k-1)!$ when $j=\tau(i)$, and with probability zero otherwise. Therefore
\[
(h_\tau)_1(x,y)=
\frac{k}{(k-1)!}\left(
\sum_{i=1}^{k}b_{i}(x)b_{\tau(i)}(y)-\frac1k
\right).
\]

Since $\sum_i b_{i}=1$, the function in parentheses corresponds to the matrix with entries $\mathbbm{1}_{\{j=\tau(i)\}}-1/k$ under the linear map
\[
(c_{ij})_{i,j\in[k]}\longmapsto
\sum_{i,j=1}^k c_{ij}b_i(x)b_j(y).
\]
These centered permutation matrices have zero row and column sums. Conversely, for any $i,j<k$, the matrix having entries $1$ at $(i,j),(k,k)$, entries $-1$ at $(i,k),(k,j)$, and zero elsewhere is a difference of two permutation matrices: match the remaining rows and columns identically in the two matrices. These $(k-1)^2$ matrices form a basis of the space with zero row and column sums, since their upper-left $(k-1)\times(k-1)$ blocks are the usual matrix units and that block determines all remaining entries. The centered permutation matrices consequently span exactly this space.

The polynomials $b_i$ are linearly independent: after division by $(1-t)^{k-1}$ on $(0,1)$, a linear relation becomes a polynomial relation in $t/(1-t)$ with distinct powers. Their products $b_i(x)b_j(y)$ are therefore linearly independent as well, so the displayed map is injective. The same holds in $L^2([0,1]^2)$, since a continuous polynomial that is zero almost everywhere is identically zero. The first-order pattern functions thus span a space of dimension $(k-1)^2$. Their covariance matrix is their $L^2$ Gram matrix, whose rank equals this dimension, proving the claim.
\end{proof}

\section{Introducing the Steering Algorithm}\label{sec:algorithm}

\subsection{Resampling Variables and First order Change}\label{sec:useful-resamplings}

We steer $U_{n}$ by resampling variables. Independently of the original sample $X_{1},\dots,X_{n}$, sample $X_{1}',\dots,X_{n}'$ independently with law $\mu$. For $D\subseteq[n]$, denote $X^{D}$ as the sample obtained by replacing $X_{i}$ with $X_{i}'$ precisely when $i\in D$. Write $U_{n}(X^{D})$ for its vector of $U$-statistics. The random algorithm outputs a family of modifications $\{D_{i}\}_{i\in[m]}$, designed so that $\{U_{n}(X^{D_i})\}_{i\in[m]}$ reaches small neighborhoods of nearby target values.

Recall from \cref{lem:first-order-scales} that $U_{n}$ can be decomposed as
\[
U_n-\mathbbm{E}[U_n]
=\frac{1}{n}\sum_{i=1}^{n}
\bigl((h_j)_1(X_i)\bigr)_{j\in[d]}+R_n.
\]
In \cref{sec:control-change}, we will show that the change induced on $R_{n}$ when resampling variables is negligible for the bounds in \cref{thm:main}, and so we design the algorithm according to the effect of resampling only on the first-order component. Explicitly, define
\[
\Gamma_i=
\bigl((h_j)_1(X_i')-(h_j)_1(X_i)\bigr)_{j\in[d]}.
\] 
The vectors $\{\Gamma_i\}_{i\in[n]}$ are independent and identically distributed, centered, and symmetric in distribution. We also have that $\Cov(\Gamma_i)=2\Cov(((h_{j})_{1}(X_{i}))_{j\in[d]})$, and so its essential linear span is $S$. Therefore, the support of $\Gamma_{i}$ contains a basis $v_1,\ldots,v_s$ of $S$. We steer the statistics by resampling variables which induce changes close to $\pm v_1,\ldots,\pm v_s$. Fix $0<\kappa<1$. For each $j\in[s]$, choose sufficiently small balls $B_j^+,B_j^-$ in $S$ around $v_j$ and $-v_j$, respectively, with one common radius, so that all $2s$ balls are disjoint and each basis coordinate of a vector in $B_j^+$ or $B_j^-$ is within $\kappa/s$ of the corresponding coordinate of $v_j$ or $-v_j$. In particular, $B_j^-=-B_j^+$. Since the centers belong to the support, $\Gamma_{i}$ takes values in each ball with positive probability, and opposite balls have equal probability by symmetry.

We say that a pair $(X_i,X_i')$ is useful in direction $j$ with positive orientation if $\Gamma_i\in B_j^+$, and with negative orientation if $\Gamma_i\in B_j^-$, and denote these sets of indices by $A_{j}^{+},A_{j}^{-}$, respectively. We only resample the variables with indices in these sets. We now move to describing how we choose the different sets of variables $\{D_{i}\}_{i\in[m]}$ that we resample.

\subsection{Resampling Prefixes of Useful Pairs}\label{sec:prefix-geometry}

Choose a uniform random ordering of $[n]$, independently of all pairs. This induces an ordering of the elements of $A_{j}^{+},A_{j}^{-}$ for each $j\in[s]$. Fix an integer $L\geq1$, representing the total number of resamplings of each type we wish to perform. We will always have $L=\Theta(\sqrt{n})$. For each vector $t=(t_1,\ldots,t_s)\in\{-L,\ldots,L\}^s$, we define a set of indices $D(t)$ whose corresponding variables we resample. If $t_j\geq0$, define $D_{j}(t_{j})$ as the first $t_{j}$ elements of $A_{j}^{+}$; if $t_j<0$, define $D_j(t_j)$ as the first $|t_j|$ elements of $A_{j}^{-}$, in both cases according to the random ordering. We say that $t$ is feasible if each required prefix exists. For feasible $t$, define $D(t)=\bigcup_{j\in[s]}D_{j}(t_{j})$ as the union of all these resampling sets. Every $t\in\{-L,\ldots,L\}^s$ is feasible whenever $L\leq\min_{j\in[s]}\min\{|A_j^+|,|A_j^-|\}$, which holds with high probability when $L=\Theta(\sqrt{n})$ (see \cref{lem:starting-conditions}). When $t$ is infeasible, we set $D(t)=\varnothing$.

Thus the random algorithm outputs $\{D(t)\}_{t\in \{-L,\dots,L\}^{s}}$. In the next section, we show that the distribution of $X^{D(t)}$ is similar enough for all $t$ when $L=\Theta(\sqrt{n})$, such that the probability of any measurable event is determined by the expected proportion of $X^{D(t)}$ satisfying that event up to a constant factor.

\section{Probability Bounds}\label{sec:probability-proof}

\subsection{Good Starting Sample and Fair Orientations}\label{sec:probability-comparison}

If every modified sample $X^{D(t)}$ had the original distribution, the expected proportion of modifications reaching a set would equal the probability that the original sample belongs to that set. Our choice of resamplings in \cref{sec:algorithm} depends on the pair values $(X_{i},X_{i}')$, so this equality need not hold. However, we only resample order $\sqrt{n}$ of the variables. For a binomial random variable with variance of order $n$, shifting a count within order $\sqrt n$ of its expectation by order $\sqrt n$ changes its point probability by at most a constant factor, and thus does not affect the estimate by too much. This analogy will hold for our case as well. To show this, we first prove that, conditional on the pair values with their useful orientations removed, the orientations of useful pairs are independent fair Bernoulli variables.

\begin{lemma}[Useful resamplings and fair orientations]\label{lem:fair-orientations}
Let $\mathcal{A}$ be the $\sigma$-algebra generated by the pairs after swapping the two entries of every negatively oriented useful pair and leaving all other pairs unchanged. Thus $\mathcal A$ retains the pair values and useful directions but excludes the orientations of useful pairs. Conditional on $\mathcal{A}$, the orientations of all useful pairs are independent fair Bernoulli variables.
\end{lemma}

\begin{proof}
Recall that $\Gamma_i$ is the change induced by replacing
$X_i$ with $X_i'$. Since $X_i$ and $X_i'$ are independent
with the same distribution, their joint law is invariant
under swapping them. Such a swap replaces $\Gamma_i$
with $-\Gamma_i$, and since $B_j^-=-B_j^+$, it preserves
whether the pair is useful and in which direction, while
reversing its orientation.

More precisely, let $T_i$ be the transformation that swaps
$X_i$ and $X_i'$ if the pair is useful, and leaves them
unchanged otherwise. This transformation preserves the
joint distribution of all pairs: the original distribution
is invariant under swapping the entries of the $i$-th
pair, and the event that this pair is useful is itself
invariant under this swap. Moreover, $T_i$ preserves
$\mathcal A$ and reverses only the orientation of the
$i$-th pair. The transformations $T_i$ commute, so we may
apply them to any collection of useful pairs without
changing either their joint distribution or $\mathcal A$.

We now show that all assignments of orientations have
equal conditional probability. Fix disjoint sets
$I_1,\ldots,I_s\subseteq[n]$, and consider the event
\[
E=\bigl\{A_j^+\cup A_j^-=I_j
\text{ for every }j\in[s]\bigr\}.
\]
This event is determined by $\mathcal A$. Write
$I=\bigcup_{j=1}^s I_j$ and $N=|I|$. For each assignment
$\varepsilon\in\{-1,1\}^I$, let $E_\varepsilon$
be the event that $E$ holds and every useful pair
$i\in I$ has orientation $\varepsilon_i$.

Given any two assignments $\varepsilon,\varepsilon'$,
swapping precisely the pairs on which they differ
maps $E_\varepsilon$ to $E_{\varepsilon'}$.
As these swaps preserve both the distribution and
$\mathcal A$, for every event $H\in\mathcal A$ we have
\[
\mathbb P(H\cap E_\varepsilon)
=
\mathbb P(H\cap E_{\varepsilon'}).
\]
The $2^N$ events $E_\varepsilon$ are disjoint and
their union is $E$. Therefore,
\[
\mathbb P(H\cap E_\varepsilon)
=2^{-N}\mathbb P(H\cap E).
\]
Since $E$ is determined by $\mathcal A$, this shows
that, conditional on $\mathcal A$, every assignment
of the $N$ useful orientations has probability
$2^{-N}$. As this holds for every possible choice
of the sets $I_1,\ldots,I_s$, the orientations of all
useful pairs are conditionally independent fair
Bernoulli variables.
\end{proof}

Therefore, conditional on $\mathcal A$, the remaining orientations are independent fair Bernoulli variables, and the number of positive orientations in each direction has a binomial distribution. The modifications described in \cref{sec:algorithm} change only the orientation of useful pairs, and thus maintain any conditioned value of $\mathcal{A}$. Thus, when the original sample is well behaved, modifications do not affect the distribution by much when conditioned on $\mathcal{A}$.

We now show that, with probability arbitrarily close to one, the sample is well behaved, meaning that the original value of $U_{n}(X)$ is sufficiently close to expectation in both the first-order component and the remainder, and that there are many useful pairs in each direction and each orientation with which to apply the algorithm. Let $N_j=|A_{j}^{+}|+|A_{j}^{-}|$ be the number of useful pairs in direction $j$, and let $K_j=|A_{j}^{+}|$ be the number with positive orientation.

\begin{lemma}[Starting Conditions]\label{lem:starting-conditions}

Assume $L=\Theta(\sqrt{n})$. For any $\epsilon>0$, there exist sufficiently large $B>0$ and sufficiently small $\nu>0$ such that, for all sufficiently large $n$,
$$
\left\|\sum_{i=1}^n\bigl((h_j)_1(X_i)\bigr)_{j\in[d]}\right\|
\leq B\sqrt n,
\qquad n\|R_n\|\leq B 
$$
$$\forall j\in[s],\qquad N_{j}\geq\nu n,\qquad |\frac{N_{j}}{2}-K_{j}|\leq B\sqrt{n}$$
all hold jointly with probability at least $1-\epsilon$. Further, there exists $c>0$ such that $N_{j}\geq\nu n$ for every $j\in[s]$ holds with probability $1-O(e^{-nc})$.
\end{lemma}

\begin{proof}
    The first two inequalities hold with probability close to one by applying Chebyshev's inequality to the bounds proven in \cref{lem:first-order-scales}. For the number of useful pairs in a given direction, each pair is useful in direction $j$ independently with positive probability not dependent on $n$, and thus this is a standard tail estimate of a binomial random variable. One can deduce it is $1-O(e^{-nc})$, for example, by Hoeffding's inequality. By \cref{lem:fair-orientations}, we have that orientations of useful pairs are independent fair Bernoulli variables for each given value of $\mathcal{A}$. Whenever $N_{j}\geq\nu$ holds, the event $|\frac{N_{j}}{2}-K_{j}|\leq B\sqrt{n}$ is the event that a binomial is at most a constant number of standard deviations from its expectation, and so Chebyshev shows this holds with probability arbitrarily close to one. A union bound on all events then shows they all hold jointly with probability at least $1-\epsilon$.
\end{proof}

\subsection{Estimating Probability by Steering Multiplicity}

For a measurable set $H\subseteq\Omega^n$, write
\[
V_H=\sum_{t\in\{-L,\ldots,L\}^s}
\mathbbm{1}_{\{t\text{ feasible}\}}
\mathbbm{1}_{\{X^{D(t)}\in H\}},
\qquad Q=(2L+1)^s.
\]
Thus $V_H/Q$ is the proportion of feasible successful modifications among all signed parameters. The following comparison bounds the probability of the original sample lying in $H$ using this proportion.

\begin{lemma}[Probability comparison for signed prefixes]\label{lem:prefix-comparison}
Fix $\nu,A,B>0$, and suppose $1\leq L\leq A\sqrt n$. Define
\[
G=\{N_j\geq\nu n\text{ for every }j\},
\qquad
C_B=\{|K_j-N_j/2|\leq B\sqrt n\text{ for every }j\}.
\]
There are constants $c,C>0$, depending only on $\nu,A,B,s$, such that, for all sufficiently large $n$ and every measurable $H\subseteq\Omega^n$,
\begin{equation}\label{eq:prefix-probability-comparison}
\frac{\mathbbm{E}[\mathbbm{1}_G\mathbbm{1}_{C_B}V_H]}{CQ}
\leq\mathbbm{P}(X\in H)
\leq\frac{\mathbbm{E}[\mathbbm{1}_GV_H]}{cQ}
+\mathbbm{P}(G^c).
\end{equation}
The expectations include the independent ordering. The constant $c$ does not depend on $B$.
\end{lemma}

\begin{proof}
Condition on $\mathcal A$. For one set of useful indices of size $N$, denote the probability of having $k$ pairs of specific orientation as 
\[
p_N(k)=2^{-N}\binom Nk,\qquad 0\leq k\leq N,
\]
and set $p_N(k)=0$ for other integers $k$. This is indeed the probability by Lemma~\ref{lem:fair-orientations}, and conditional on count $k$, the set of indices with that orientation is a uniform $k$-subset of the class.

Fix a count vector $k=(k_1,\ldots,k_s)$, and consider
a parameter $t\in\{-L,\ldots,L\}^s$. Conditional on
$\mathcal A$ and the counts $(K_j)_{j\in[s]}=k$, the
set $A_j^+$ is a uniformly random $k_j$-subset of
the $N_j$ useful pairs in direction $j$. If $t_j\geq 0$,
the modification reverses the orientations of the first
$t_j$ positive pairs in the independent random ordering.
Thus, it removes a uniformly random $t_j$-subset of
$A_j^+$. If $t_j<0$, it similarly adds a uniformly
random $|t_j|$-subset of $A_j^-$ to $A_j^+$.

In either case, the resulting set of positively oriented
pairs is uniform among all subsets of size $k_j-t_j$.
Indeed, this follows from the uniformity of the original
set and the random ordering: every permutation of the
$N_j$ indices preserves their joint distribution and
commutes with the modification, so every resulting subset
of the prescribed size has the same probability.
The relative orderings and orientations of different
classes are independent, so the resulting sets are
jointly uniform over all orientation assignments with
count vector $k-t$.

Note also that the modification is feasible precisely
when $0\leq k_j-t_j\leq N_j$ for every $j\in[s]$.
Therefore, if we define
\[
f(k)=\mathbbm{P}\bigl(
X\in H\mid\mathcal A,(K_j)_{j\in[s]}=k
\bigr),
\]
then, averaging over the initial orientations and the
independent ordering, we obtain
\[
\begin{aligned}
&\mathbbm{E}\bigl[
\mathbbm{1}_{\{t\text{ feasible}\}}
\mathbbm{1}_{\{X^{D(t)}\in H\}}
\mid\mathcal A,(K_j)_{j\in[s]}=k
\bigr]\\
&\qquad=
\mathbbm{1}_{\{0\leq k_j-t_j\leq N_j\ \forall j\}}
f(k-t).
\end{aligned}
\]
This is the key observation: although the success of
a modification may depend on the specific indices of
each orientation, after averaging over the orientations
and ordering, its probability depends on the initial
orientation counts only through $k-t$.

Put $p(k)=\prod_{j=1}^s p_{N_j}(k_j)$. All sums over
$k$ below range over $\prod_{j=1}^s\{0,\ldots,N_j\}$.
By the law of total probability,
\[
\mathbbm{P}(X\in H\mid\mathcal A)
=\sum_k f(k)p(k).
\]
Similarly, summing the preceding identity over all
parameters $t$ gives
$$\mathbbm{E}[V_H\mid\mathcal A]
=\sum_k p(k)
  \sum_{\substack{t\in\{-L,\ldots,L\}^s\\
  0\leq k_j-t_j\leq N_j\ \forall j}}
  f(k-t)
=\sum_k f(k)\prod_{j=1}^s
  \left(\sum_{u=-L}^L p_{N_j}(k_j+u)\right).$$
Here the second equality follows by changing variables
from $k-t$ to $k$, and using the product structure of
$p(k)$.
 To make sure shifting the number of pairs of a specific orientation doesn't change the probability by much, we restrict the binomial count of positive orientations to be close to expectation, such that a small shift doesn't change the probability by much. Thus, we look at
\[
\mathbbm{E}[\mathbbm{1}_{C_B}V_H\mid\mathcal A]
=\sum_k f(k)\prod_{j=1}^s
\left(\sum_{u=-L}^L p_{N_j}(k_j+u)
\mathbbm{1}_{\{|k_j+u-N_j/2|\leq B\sqrt n\}}\right).
\]

These nearby central binomial probabilities are indeed comparable. Suppose $N\geq\nu n$ and $|k-N/2|\leq D\sqrt n$, where $D$ is fixed. A direct and standard computation for binomial distributions shows that there exist $c_{0},C_{0}$, depending only on $D$ such that
\[
\sum_{u=-L}^L p_N(k+u)\geq c_0(2L+1)p_N(k),\qquad 
\sum_{u=-L}^L p_N(k+u)
\mathbbm{1}_{\{|k+u-N/2|\leq B\sqrt n\}}
\leq C_0(2L+1)p_N(k).
\]

Apply these two estimates to each class on $G$, which is determined by $\mathcal A$, and then average over $\mathcal A$. We obtain
\[
\mathbbm{E}[\mathbbm{1}_GV_H]
\geq c_0^s Q\,\mathbbm{P}(X\in H,G),
\qquad
\mathbbm{E}[\mathbbm{1}_G\mathbbm{1}_{C_B}V_H]
\leq C_0^s Q\,\mathbbm{P}(X\in H,G).
\]
Since $\mathbbm{P}(X\in H)\leq\mathbbm{P}(X\in H,G)+\mathbbm{P}(G^c)$, these are the claimed bounds.
\end{proof}

The restriction to typical starting counts is needed only for the lower comparison. Without it, shifts from nearer the middle can substantially increase the probability of a rare extreme count.

Now that we have established that the number of modifications satisfying an event bounds the probability of the original sample satisfying it, we move to analyze $V_{H}$ in the case of values of $U_{n}$. 

\section{Analyzing the Change Induced by Steering}\label{sec:control-change}

We handle the change induced on the first-order component of $U_{n}$ and the remainder of $U_{n}$ separately, showing first that the change to the first-order component is controlled, and then that the change to the remainder is small with high probability. The estimates of the first order apply whenever both lists in every direction have length at least $L$. The moment estimates for the remainder are unconditional.

\subsection{Control of The First-Order Change}\label{sec:first-order-control}

We now show that these resamplings induce first-order changes which cover a central region up to a bounded error, and separate different parameters. For convenience, we identify $S$ with its coordinates in the basis $v_1,\ldots,v_s$. Let $F(t)$ be the basis-coordinate vector of $\sum_{i\in D(t)}\Gamma_i$ at integer parameters. When both lists in every direction have length at least $L$, extend this map to $[-L,L]^s$ by linearly interpolating each signed coordinate path and adding the resulting paths.

\begin{lemma}[Covering and separation by signed prefixes]\label{lem:prefix-geometry}
For all $t,u\in[-L,L]^s$,
\[
\|F(t)-F(u)-(t-u)\|_\infty
\leq\kappa\|t-u\|_\infty,
\qquad
\|F(t)-F(u)\|_\infty
\geq(1-\kappa)\|t-u\|_\infty.
\]
Consequently, every $q\in\mathbb R^s$ with $\|q\|_\infty\leq(1-\kappa)L/2$ is within $(1+\kappa)/2$ in supremum norm of $F(t)$ for some integer $t\in\{-L,\ldots,L\}^s$. Furthermore, for any $q\in\mathbb R^s$ and $B\geq0$, the number of integer parameters satisfying $\|F(t)-q\|_\infty\leq B$ is at most $K(1+B)^s$, where $K$ depends only on $s$ and $\kappa$.
\end{lemma}

\begin{proof}
Let $e_j$ denote the $j$th coordinate vector. Along the positive part of the $j$th path, each unit increment of the parameter adds a response whose coordinates differ from $e_j$ by at most $\kappa/s$. Along its negative part, increasing the parameter by one removes a response close to $-e_j$, so the same bound holds. Linear interpolation preserves this bound on each segment. Splitting an interval at zero and at integer points thus gives, for all $a,b\in[-L,L]$, including opposite signs,
\[
\|F(ae_j)-F(be_j)-(a-b)e_j\|_\infty
\leq\frac{\kappa}{s}|a-b|.
\]
The map is the sum of these $s$ coordinate paths. Summing the last inequality proves the first displayed bound, and the reverse triangle inequality proves the second.

For covering, fix $q$ in the stated central cube and consider, on $[-L/2,L/2]^s$, the map
\[
t\longmapsto q-(F(t)-t).
\]
Since $F(0)=0$, the first bound shows that this map takes the cube into itself and has Lipschitz constant at most $\kappa<1$. Its unique fixed point satisfies $F(t)=q$. Round each coordinate to a nearest integer, obtaining $u$. Then $\|u-t\|_\infty\leq1/2$, and $u\in[-L,L]^s$ because $L\geq1$. The first bound also gives $\|F(u)-F(t)\|_\infty\leq(1+\kappa)/2$, as required.

Finally, if the set of integer parameters in question is nonempty, fix one of them, $u$. Every other such $t$ satisfies
\[
(1-\kappa)\|t-u\|_\infty
\leq\|F(t)-F(u)\|_\infty\leq2B.
\]
There are at most $(1+4B/(1-\kappa))^s$ integer points in this box.
\end{proof}

\subsection{Uniform Control of The Remainder}\label{sec:resampling-responses}

The first-order decomposition gives, simultaneously for all $D\subseteq[n]$,
\[
U_n(X^D)-U_n(X)
=\frac1n\sum_{i\in D}\Gamma_i+R_n(X^D)-R_n(X).
\]
In \cref{lem:first-order-scales}, we showed that $R_{n}$ has small variance. We strengthen this by showing that the change in $R_{n}$ by resampling is small simultaneously for many choices of modification sets $D$. We obtain different bounds when the kernels are bounded compared to when they are only of finite variance, which is what brings about the different upper bounds in \cref{thm:main}.

\begin{lemma}[Uniform resampling error]\label{lem:uniform-resampling-error}
Let $R\subseteq[n]$ have cardinality $m\leq n/2$, and be independent of all variables. For $n\geq4\max_j k_j$, there is a constant $K$, depending only on the fixed kernels and their orders, such that
\[
\left\|
\max_{D\subseteq R}
\left\|
U_n(X^D)-U_n(X)
-\frac{1}{n}\sum_{i\in D}\Gamma_i
\right\|
\right\|_{L^2}
\leq
K\left(\frac{m}{n^{3/2}}+\frac{m^2}{n^2}\right).
\]
For every fixed even integer $p\geq2$, the same bound holds in $L^p$ if the first-order-degenerate remainder of each kernel belongs to $L^p$, with the constant depending on $p$. In particular, if $m=O(\sqrt n)$, the uniform error is $O_{L^2}(n^{-1})$; for bounded kernels it is $O_{L^p}(n^{-1})$ for every fixed even $p$.
\end{lemma}

\begin{proof}
Condition on $R$. All estimates below are uniform in its choice. It suffices to prove the claim for one kernel $h$ of order $k$, and then apply the triangle inequality to the coordinatewise bounds. Fix an even integer $p\geq2$ for which the remainder kernel $g=g_h$ belongs to $L^p$. If $k=1$, then $g=0$, so assume $k\geq2$. Writing the remainder as $U_{n}(g)$, the change in the remainder which we wish to bound is of the form
$$U_n(g;X^D)-U_n(g;X)=\binom{n}{k}^{-1}\sum_{I\in\binom{[n]}{k}}\bigl(g(X_{I};X^{D})-g(X_{I})\bigr).$$
We partition the sets $I\in\binom{[n]}{k}$ according to $|I\cap R|$. When $|I\cap R|=0$, the difference by resampling is zero. When $|I\cap R|=1$ we have jointly for all $D\subseteq R$ that
$$\left|\sum_{i\in R}\sum_{I\in \binom{[n]\setminus R}{k-1}}\bigl(g(X_{I\cup\{i\}};X^{D})-g(X_{I\cup\{i\}})\bigr)\right|\leq\sum_{i\in R}\left|\sum_{I\in\binom{[n]\setminus R}{k-1}}\bigl(g(X_{i}',X_{I})-g(X_{i},X_{I})\bigr)\right|.$$
Therefore, this bound holds for the maximum over all $D\subseteq R$ as well. For $|I\cap R|\geq 2$, even when we do not resample all variables in the intersection the corresponding sets may still induce a change. Therefore, we only have the crude bound of
$$\left|\sum_{\substack{I\in\binom{[n]}{k}\\|I\cap R|\geq 2}}\bigl(g(X_{I};X^{D})-g(X_{I})\bigr)\right|\leq \sum_{\substack{I\in\binom{[n]}{k}\\|I\cap R|\geq 2}}\sum_{J\subseteq I\cap R}|g(X_I;X^J)-g(X_{I})|.$$
Here we sum over all possibilities of resampled subsets for each $I$, and take absolute values of all summands. Though much coarser than the bound for sets with intersection of size one, this will be sufficient for the bound as there are fewer such sets when $R$ is small. Combining these two bounds, which apply jointly over all $D\subseteq R$, we get that $\max_{D\subseteq R}|U_{n}(g;X^{D})-U_{n}(g;X)|$ is bounded from above by
$$\binom{n}{k}^{-1}\left(\sum_{i\in R}\left|\sum_{I\in\binom{[n]\setminus R}{k-1}}\bigl(g(X_{i}',X_{I})-g(X_{i},X_{I})\bigr)\right|+\sum_{\substack{I\in\binom{[n]}{k}\\|I\cap R|\geq 2}}\sum_{J\subseteq I\cap R}|g(X_I;X^J)-g(X_{I})|\right).$$
Therefore, it is sufficient to bound the $L^{p}$ norm of this expression. We start with the case $|I\cap R|\geq 2$. For each such set, the triangle inequality and Minkowski's inequality give us
$$\left\|\sum_{J\subseteq I\cap R}|g(X_I;X^J)-g(X_{I})|\right\|_{L^{p}}\leq \sum_{J\subseteq I\cap R}\left(\|g(X_I;X^J)\|_{L^{p}}+\|g(X_{I})\|_{L^{p}}\right)\leq 2^{k+1}\|g\|_{L^{p}}$$
as all kernel evaluations have the same distribution and thus the same $L^p$ norm. Since there are at most $\binom{m}{2}\binom{n}{k-2}$ such sets, Minkowski's inequality shows that the $L^p$ norm of the second sum, after multiplication by $\binom{n}{k}^{-1}$, is at most $K\frac{m^{2}}{n^{2}}\|g\|_{L^p}$, where $K$ depends only on $k$.

We move to the case $|I\cap R|=1$. We fix $i\in R$, and look at
$$\mathbbm{E}\left[\left(\sum_{I\in\binom{[n]\setminus R}{k-1}}\bigl(g(X_{i}',X_{I})-g(X_{i},X_{I})\bigr)\right)^{p}\right]=\sum_{I_{1},\dots,I_{p}\in\binom{[n]\setminus R}{k-1}}\mathbbm{E}\left[\prod_{j\in[p]}\bigl(g(X_{i}',X_{I_{j}})-g(X_{i},X_{I_{j}})\bigr)\right].$$
Here, there are too many summands for naively bounding each one to be helpful, and so we first show that most $p$-tuples give zero expectation. Indeed, suppose one $I_j$ is disjoint from all the others. Conditional on $X_i,X_i'$ and all variables indexed by the other sets, its factor has mean zero by degeneracy of $g$, while all other factors are fixed. Therefore, it is sufficient to count the number of $p$-tuples in which for every $j\in[p]$ we have $I_j\cap\bigl(\bigcup_{l\in[p]\setminus\{j\}}I_l\bigr)\neq\emptyset$. Equivalently, there are no isolated vertices in the intersection graph of $I_1,\dots,I_p$, in which each set is a vertex and two sets are connected by an edge if they intersect. If a connected component of this graph has $r$ sets, we can adjoin them consecutively so that each set after the first intersects at least one of the previous ones. Each set starting from the second therefore contributes at most $k-2$ new elements, so the component has at most $(k-1)+(r-1)(k-2)=(k-2)r+1$ distinct elements. As there are no isolated vertices, each connected component has size at least $2$ and there are at most $p/2$ connected components. Therefore, denoting their sizes by $r_1,\dots,r_a$, the number of distinct elements in such a tuple is at most
$$\sum_{l=1}^{a}\bigl((k-2)r_l+1\bigr)=(k-2)p+a\leq(k-2)p+\frac p2=(k-1)p-\frac p2.$$
For each choice of at most $(k-1)p-\frac{p}{2}$ indices, there is a number of options, bounded in terms of $k,p$ alone, for a tuple $I_{1},\dots,I_{p}$ whose union is contained in those indices. Therefore, the number of potentially nonzero summands is $O(n^{(k-1)p-\frac{p}{2}})$. By H\"older's inequality, the absolute value of each product expectation is at most $(2\|g\|_{L^p})^p$. Taking the $p$-th root of the moment therefore shows that, for each $i\in R$, the inner sum has $L^p$ norm $O(n^{k-\frac32}\|g\|_{L^p})$. Applying Minkowski's inequality to the entire first sum, we have
$$\left\|\binom{n}{k}^{-1}\sum_{i\in R}\left|\sum_{I\in\binom{[n]\setminus R}{k-1}}\bigl(g(X_{i}',X_{I})-g(X_{i},X_{I})\bigr)\right|\right\|_{L^p}\leq K\frac{m}{n^{\frac{3}{2}}}\|g\|_{L^p}$$
for a sufficiently large $K$ depending only on $k,p$. Combining the two bounds gives the desired result. 
\end{proof}

\section{Proving the Bounds in \cref{thm:main}}\label{sec:main-proof}

\subsection{The finite-variance lower bound}\label{sec:probability-lower}

If $s=0$, the only permitted target is $z=0$, and Lemma~\ref{lem:first-order-scales} gives
\[
\mathbbm{P}\left(\|U_n-\mathbbm{E}[U_n]\|>C/n\right)\leq K/C^2.
\]
Choosing $C$ sufficiently large proves the lower bound. Suppose henceforth that $s>0$, and fix $M$. Note first that by \cref{lem:starting-conditions}, the original sample starts sufficiently close to all permitted targets, and $G,C_{B}$ are satisfied for large enough $B$ and small enough $\nu$. Thus, choosing $L$ to be large enough of order $\sqrt{n}$, we have by \cref{lem:prefix-geometry}, that for every $z\in S$ with $\|z\|\leq M/\sqrt n$, the vector
\[
nz-\sum_{i=1}^n\bigl((h_j)_1(X_i)\bigr)_{j\in[d]}
\]
belongs to the covering region, and there exists some prefix $t$ such that the change in first order component of $U_{n}$ steers its value to within order $n^{-1}$ of the prescribed value. By \cref{lem:uniform-resampling-error}, the maximal shift in the remainder over resampling of all subsets of the first $A\sqrt{n}$ variables of the random ordering has $L_{2}$ norm of order at most $n^{-1}$. Assuming that the first $A\sqrt{n}$ indices in the random ordering contain $L$ useful pairs of each direction and orientation, which occurs with probability arbitrarily close to one for large enough $A$, we have by Chebyshev that the change to the remainder for the prefix $t$ providing the desired first order affect is of order $n^{-1}$. Thus, applying the orthogonal decomposition presented in \cref{lem:first-order-scales}, under these assumptions there exists such a prefix $t$ and constant $C$ satisfying
\[
\left\|U_n(X^{D(t)})-\mathbbm{E}[U_n]-z\right\|
\leq\frac{C}{n}.
\]
note that we bounded the change to the remainder throughout all $\mathbb{R}^{d}$, and so this also controls the coordinates in $S^\bot$. We look at the event
\[
H_z=\{x\in\Omega^n:\|U_n(x)-\mathbbm{E}[U_n]-z\|\leq C/n\}.
\]
The preceding discussion proves that $V_{H_{z}}\geq 1$ given $G,C_{B}$, and hence
\[
\mathbbm{E}[\mathbbm{1}_G\mathbbm{1}_{C_{B_0}}V_{H_z}]\geq \frac{1}{2}.
\]
Applying the lower comparison in Lemma~\ref{lem:prefix-comparison}, we obtain
\[
\mathbbm{P}\left(\|U_n-\mathbbm{E}[U_n]-z\|\leq C/n\right)
\geq\frac{1}{2C(2L+1)^s}\geq c n^{-s/2}.
\]
The good event and all constants were chosen independently of $z$, proving the required uniformity.

\subsection{Upper bounds from separation and multiplicity}\label{sec:probability-upper}

Throughout this section, we assume both the high probability events of \cref{lem:starting-conditions}, that $L,A=\Theta(\sqrt{n})$ are sufficiently large, and that there exists $L$ useful pairs of each direction and orientation in the first $A$ pairs of the random ordering.

For the upper bound, we use the other side of the comparison in \cref{lem:prefix-comparison}. The separation result from \cref{lem:prefix-geometry} limits the number of prefixes $t$ steering the values of the first order component of $U_{n}$ to a prescribed small ball to a constant. Therefore, the multiplicity can only grow when the change to the remainder is large. Denote the scaled uniform resampling change in remainder as
$$Z_{n}=n\max_{t\in \{-L,\dots,L\}^{s}}
\left\|
U_n(X^{D(t)})-U_n(X)
-\frac{1}{n}\sum_{i\in D}\Gamma_i
\right\|.$$ 
We will show that the increase in multiplicity due to shift of the remainder is bounded by $(1+Z_{n})^{s}$, which will allow to apply the moment bounds of \cref{lem:uniform-resampling-error} to achieve an upper bound.

\begin{lemma}[Multiplicity of signed-prefix endpoints]\label{lem:prefix-multiplicity}
Whenever $R$ contains at least $L$ useful pairs of each direction and orientation, for every $C>0$ and every $a\in S$ we have
\[
\left|\left\{t\in\{-L,\ldots,L\}^s:
\|P_SU_n(X^{D(t)})-a\|\leq C/n\right\}\right|
\leq\min\left\{(2L+1)^s,K_C(1+Z_n)^s\right\}.
\]
The constant depends only on $C$ and the fixed response classes.
\end{lemma}

\begin{proof}
If $t,u$ are two prefixes which steer $U_{n}$ to the prescribed set of values, the orthogonal decomposition gives us that the difference in the change they induced on the first order component is at most
\[
\left\|\sum_{i\in D(t)}\Gamma_i-
\sum_{i\in D(u)}\Gamma_i\right\|\leq2(C+Z_n).
\]
By Lemma~\ref{lem:prefix-geometry}, this bounds $\|t-u\|_\infty$ by a constant times $C+Z_n$. Therefore, whenever there is a prefix steering to the desired set of values, the collection of all prefixes steering to that set of values lies in the box of diameter of order $C+Z_{n}$ in each coordinate. Consequently, the number of such prefixes is at most $K_C(1+Z_n)^s$.
\end{proof}

For $z\in S$, we apply the upper comparison in Lemma~\ref{lem:prefix-comparison} to
\[
H_z=\{x\in\Omega^n:\|P_SU_n(x)-a\|\leq C/n\}.
\]
 By lemma~\ref{lem:prefix-multiplicity}, $V_{H_a}$ is bounded by $\min\{Q,K_C(1+Z_n)^s\}$. By \cref{lem:starting-conditions}, we have $\mathbbm{P}(G^c)=O(e^{-cn})$. Consequently,
\[
\mathbbm{P}(X\in H_z)
\leq\frac{1}{cQ}\mathbbm{E}\left[\min\{Q,K_C(1+Z_n)^s\}\right]
+O(e^{-c'\sqrt n}).
\]
Since $Q\asymp L^s$, this provides
\begin{equation}\label{eq:projected-multiplicity-bound}
\sup_{a\in S}\mathbbm{P}\left(\|P_SU_n-a\|\leq C/n\right)
\leq K L^{-s}\mathbbm{E}\left[\min\{L^s,K_C(1+Z_n)^s\}\right]
+O(e^{-c\sqrt n}),
\qquad L\asymp\sqrt n.
\end{equation}

For bounded kernels, Lemma~\ref{lem:uniform-resampling-error} bounds the moments of $Z_n$ of every fixed even order. Choose such an order at least $s$. Then $\mathbbm{E}[(1+Z_n)^s]=O(1)$, and \eqref{eq:projected-multiplicity-bound} gives $O(L^{-s})=O(n^{-s/2})$. The same argument applies whenever the first-order-degenerate remainder kernels have a finite even moment of order at least $\max\{2,s\}$; the first-order projections themselves need only have finite second moments.

Under finite variance alone, $\mathbbm{E}[Z_n^2]=O(1)$ gives the same bound for $s=1,2$. When $s>2$, the elementary inequality
\[
\min\{L^s,b^s\}\leq L^{s-2}b^2,\qquad b\geq0,
\]
first gives $O(L^{-2})=O(n^{-1})$.

To obtain the little-$o$ improvement, we approximate the remainder kernels by bounded ones. Recall that by \cref{lem:first-order-scales}, the remainder $R_n$ is $(U_n(g_j))_{j\in[d]}$ for centered and degenerate kernels of the form
\[
g_j(x_1,\ldots,x_{k_j})=
h_j(x_1,\ldots,x_{k_j})-\mathbbm{E}[h_j(X_1,\ldots,X_{k_j})]
-\frac1{k_j}\sum_{a=1}^{k_j}(h_j)_1(x_a).
\]
Each $g_j$ is centered and has vanishing first-order projection. For $B>0$, put
\[
f_{j,B}=\max\{-B,\min\{g_j,B\}\},
\qquad
g_{j,B}(x_1,\ldots,x_{k_j})
=f_{j,B}(x_1,\ldots,x_{k_j})-\mathbbm{E}[f_{j,B}]
-\frac1{k_j}\sum_{a=1}^{k_j}(f_{j,B})_1(x_a).
\]
Thus $g_{j,B}$ is bounded, symmetric, centered and first-order-degenerate. We check that $g_{j,B}\to g_j$ in $L^2$. Writing $q_{j,B}=g_j-f_{j,B}$ and evaluating kernels at $(X_1,\ldots,X_{k_j})$, the centering and degeneracy of $g_j$ give
\[
g_j-g_{j,B}
=q_{j,B}-\sum_{a=1}^{k_j}\mathbbm{E}[q_{j,B}\mid X_a]
+(k_j-1)\mathbbm{E}[q_{j,B}].
\]
Conditional Jensen's inequality and the triangle inequality therefore imply
\[
\|g_j-g_{j,B}\|_{L^2}
\leq2k_j\|q_{j,B}\|_{L^2}\longrightarrow0.
\]
Set $\delta_B=\sum_{j=1}^d\|g_j-g_{j,B}\|_{L^2}$, so $\delta_B\to0$.

Let $Z_{n,B}$ and $Z_{n,B}^{\mathrm{tail}}$ be $n$ times the maximum, over all $D\subseteq R$, of the norm of the change in the vectors of $U$-statistics with kernels $(g_{j,B})_{j\in[d]}$ and $(g_j-g_{j,B})_{j\in[d]}$, respectively. The first-order terms in the original kernels cancel exactly from the resampling error, so
\[
Z_n\leq Z_{n,B}+Z_{n,B}^{\mathrm{tail}}.
\]
All these remainder kernels have zero first-order response. We may therefore apply Lemma~\ref{lem:uniform-resampling-error} directly to their changes. The linear dependence on the kernel norms in its proof gives
\[
\|Z_{n,B}^{\mathrm{tail}}\|_{L^2}
\leq K'\left(\frac m{\sqrt n}+\frac{m^2}{n}\right)\delta_B
\leq K\delta_B,
\]
uniformly for sufficiently large $n$. In particular, $\mathbbm{E}[(Z_{n,B}^{\mathrm{tail}})^2]\leq K^2\delta_B^2$. For each fixed $B$, all fixed even moments of $Z_{n,B}$ are uniformly bounded.

Now apply \eqref{eq:projected-multiplicity-bound} ,using $\min\{L^s,(Z_{n,B}^{\mathrm{tail}})^s\}\leq L^{s-2}(Z_{n,B}^{\mathrm{tail}})^2$. We obtain
\[
\sup_{z\in S}\mathbbm{P}\left(\|P_SU_n-z\|\leq C/n\right)
\leq K_B n^{-s/2}+K\delta_B^2 n^{-1}+O(e^{-c\sqrt n}),
\]
where $K$ is independent of $B$. For each fixed $B$, multiply by $n$ and take the upper limit as $n\to\infty$. Since $s>2$, the result is at most $K\delta_B^2$. Letting $B\to\infty$ proves the uniform $o(n^{-1})$ bound. Only the estimate of the remainder error was changed in this argument; the response law, modifications and first-order fluctuation space remained fixed. This proves all the probability bounds in \cref{thm:main}. Its sharpness assertion is proved in \cref{sec:sharpness}.

\section{Proof of \cref{thm:permutations}}\label{sec:permutation-proof}

Apply the representation of \cref{sec:permutation-prelim} to a uniform permutation $\Pi_n$. For its vector of normalized $k$-pattern counts, with $k!$ coordinates, denote
\[
D_k(\Pi_n)=\binom nk
\max_{\tau\in S_k}
\left|U_n(h_\tau)-\mathbbm{E}[U_n(h_\tau)]\right|.
\]
The kernels are bounded and, by \cref{lem:permutation-rank}, their first-order fluctuation rank is $(k-1)^2$.

Let $C_0$ be a radius constant in the lower bound of \cref{thm:main} for this vector and $M=0$. On the event $\|U_n-\mathbbm{E}[U_n]\|\leq C_0/n$, every coordinate has absolute deviation at most $C_0/n$, and hence
\[
D_k(\Pi_n)\leq\frac{C_0}{n}\binom nk
\leq\frac{C_0}{k!}n^{k-1}.
\]
Choose $C_k\geq C_0/k!$. For every fixed $C\geq C_k$, the probability that $D_k(\Pi_n)\leq Cn^{k-1}$ is therefore at least a positive constant times $n^{-(k-1)^2/2}$ for all sufficiently large $n$.

Conversely, the event $D_k(\Pi_n)\leq Cn^{k-1}$ implies
\[
\|U_n-\mathbbm{E}[U_n]\|
\leq\sqrt{k!}\,C\frac{n^{k-1}}{\binom nk}
\leq\frac{2\sqrt{k!}\,Ck!}{n}
\]
for all sufficiently large $n$, since $n^{k-1}/\binom nk\sim k!/n$. The bounded-kernel upper bound in \cref{thm:main} applies with this fixed radius. In particular,
\[
\mathbbm{P}\bigl(D_k(\Pi_n)\leq C_kn^{k-1}\bigr)
=\Theta_k\bigl(n^{-(k-1)^2/2}\bigr),
\]
and for every fixed $C\geq C_k$ the same estimate holds with constants allowed to depend on $k,C$. Since each permutation has probability $1/n!$, multiplication by $n!$ gives the enumeration assertion of \cref{thm:permutations}.

The positive lower probability also gives, for every sufficiently large $n$, a permutation with discrepancy at most $C_kn^{k-1}$. For $k\geq4$, Beniamini, Lavee and Linial~\cite{BLL25} prove the complementary bound $\min_{\pi\in S_n}D_k(\pi)=\Omega_k(n^{k-1})$. Together these bounds yield $\min_{\pi\in S_n}D_k(\pi)=\Theta_k(n^{k-1})$, completing the proof.

\section{Proving the Sharpness of \cref{thm:main}}\label{sec:sharpness}

We prove that the finite-variance upper bound of
\cref{sec:probability-upper} admits no universal rate of improvement.
For every integer $s>2$, every positive sequence $\epsilon_n\to0$,
and every fixed $C>0$, we construct a fixed symmetric square-integrable
kernel with values in $\mathbb R^s$, whose first-order fluctuation space
has dimension $s$, such that
\[
\limsup_{n\to\infty}\frac{n}{\epsilon_n}
\mathbbm{P}\left(
\|U_n-\mathbbm{E}[U_n]\|\leq\frac{C}{n}
\right)=\infty.
\]

The example uses a scalar $Z\in\mathbb{R}$ and a bounded vector $Y\in\mathbb R^s$
satisfying $\mathbbm{E}[Z]=0$, $\mathbbm{E}[Y]=0$, and $ZY=Y$ almost surely. This algebraic identity implies that whenever $Z\neq 1$ we have $Y=0$, meaning that a single scalar constraint collapses all $S$ directions. This is what we will use to prove concentration. For our $U$-statistics, we use the estimator of the covariance between $Y,Z$:
\[
h\bigl((z,y),(z',y')\bigr)=\frac12(z-z')(y-y').
\]
Applying $\mathbbm{E}[Z]=\mathbbm{E}[Y]=0$ and the identity $ZY=Y$, a direct computation gives $\mathbbm{E}[U_n]=0$ and
\begin{equation}\label{eq:sharpness-factorization}
nU_n=\frac{n-\sum_{i=1}^n Z_i}{n-1}\sum_{i=1}^n Y_i.
\end{equation}
Thus, if $\sum_i Z_i=n$, the statistic becomes zero. We give $Z$ a heavy tail so that one
observation can equal $n$ with sufficiently large probability. Note that this requires that $Z$ is not bounded, as shown to be necessary in \cref{thm:main}.
The remaining centered scalar sum and the vector sum then fluctuate
at order $\sqrt n$, leaving $U_n$ at order $n^{-1}$.

Fix $s,C$, and $(\epsilon_n)$ as above. Set $\delta=C/2$ and
$a_j=2^{-j}$, and choose strictly increasing integers
$n_j\geq2^{j^2+1}$ such that $\epsilon_{n_j}\leq a_j/j^2$.
Let $W$ have the symmetric law
\[
\mathbbm{P}(W=n_j+1)=\mathbbm{P}(W=-(n_j+1))
=\frac{a_j}{2(n_j+1)^2},\qquad j\geq1,
\]
with the remaining probability at zero. Let $\xi$ be independent of
$W$ and uniform on $\{\pm e_1,\ldots,\pm e_s\}$. With an independent
fair choice, define
\[
(Z,Y)=
\begin{cases}
(1,\delta\xi),&\text{with probability }1/2,\\
(-1+W,0),&\text{with probability }1/2.
\end{cases}
\]
The tail masses were chosen so that
$\mathbbm{E}[W]=0$ and $\mathbbm{E}[W^2]=\sum_j a_j=1$.
Consequently,
\[
\mathbbm{E}[Z]=0,\qquad \mathbbm{E}[Y]=0,\qquad ZY=Y,
\qquad \mathbbm{E}[Z^2]=\frac32.
\]
Since $Y$ is bounded, $h$ is square-integrable. Moreover, a direct computation gives
\[
(h)_1(Z,Y)=Y,\qquad
\Cov(Y)=\frac{\delta^2}{2s}I_s,
\]
so the first-order fluctuation space has dimension exactly $s$.

Now take $n=n_j$ and write $A_i=\{Z_i=n\}$. Then
\[
p_j:=\mathbbm{P}(A_i)=\frac{a_j}{4(n_j+1)^2}.
\]
On $A_i$, we have $Y_i=0$, and
\eqref{eq:sharpness-factorization} becomes
\[
nU_n=-\frac{T}{n-1}\sum_{\ell\ne i}Y_\ell,
\qquad T=\sum_{\ell\ne i}Z_\ell.
\]
Conditional on all scalar observations, the vectors $Y_\ell$ are
independent and centered, with norms at most $\delta$. Hence
\[
\mathbbm{E}\left[
\left\|\sum_{\ell\ne i}Y_\ell\right\|^2
\,\middle|\, Z_1,\ldots,Z_n\right]\leq\delta^2(n-1).
\]
The other observations retain their original law conditional on $A_i$,
so $\mathbbm{E}[T^2\mid A_i]=3(n-1)/2$. Therefore
\[
\mathbbm{E}[\|nU_n\|^2\mid A_i]
\leq\frac32\delta^2=\frac{3C^2}{8},
\qquad
\mathbbm{P}(\|U_n\|\leq C/n\mid A_i)\geq\frac58.
\]
The events $A_i$ are independent. Applying the first two terms of
inclusion--exclusion to their intersections with
$\{\|U_n\|\leq C/n\}$ gives
\[
\mathbbm{P}(\|U_{n_j}\|\leq C/n_j)
\geq\frac58 n_jp_j-\binom{n_j}{2}p_j^2
\geq\frac12 n_jp_j
\geq\frac{a_j}{32n_j},
\]
where we used $n_jp_j\leq1/16$. It follows that
\[
\frac{n_j}{\epsilon_{n_j}}
\mathbbm{P}(\|U_{n_j}-\mathbbm{E}[U_{n_j}]\|\leq C/n_j)
\geq\frac{j^2}{32}\longrightarrow\infty.
\]
All sample sizes use the same distribution and kernel.

This example also illustrates the source of multiplicity in
\cref{lem:prefix-multiplicity}. Choose the response basis
$v_j=2\delta e_j$ and sufficiently small response balls. The useful
pairs are precisely the ordered pairs
$\bigl((1,-\delta e_j),(1,\delta e_j)\bigr)$ and their reversals,
so steering preserves every $Z_i$. For every feasible prefix,
\[
\begin{aligned}
n\bigl(U_n(X^{D(t)})-U_n(X)\bigr)
&=\frac{n-\sum_i Z_i}{n-1}\,2\delta t,\\
n\bigl(R_n(X^{D(t)})-R_n(X)\bigr)
&=\frac{1-\sum_i Z_i}{n-1}\,2\delta t.
\end{aligned}
\]
When $\sum_i Z_i=n+O(\sqrt n)$, all prefixes with
$\|t\|_\infty\leq L=\Theta(\sqrt n)$ produce endpoints within
$O(n^{-1})$ of one another, while the scaled change in the remainder
can be of order $\sqrt n$. The large multiplicity comes from this
exceptional cancellation of the first-order responses by the remainder.

\bibliographystyle{alpha}
\bibliography{local-steering}
\end{document}